\documentstyle{amsppt}
\magnification\magstep1
\NoRunningHeads
\baselineskip=15pt
\input btxmac.tex               
\bibliographystyle{plain}
\loadbold
\input diagrams



\def\ideal#1.{I_{#1}}
\def\ring#1.{\Cal O_{#1}}
\def\proj#1.{\Bbb P(#1)}
\def\pr #1.{\Bbb P^{#1}}
\def\af #1.{\Bbb A^{#1}}
\def\Hz #1.{\Bbb F_{#1}}
\def\Hbz #1.{\overline{\Bbb F}_{#1}}
\def\fb#1.{\underset #1 \to \times}
\def\ten#1.{\underset #1 \to \otimes}
\def\res#1.{\underset {\ \ring #1.} \to \otimes}
\def\au#1.{\operatorname {Aut}\,(#1)}
\def\deg#1.{\operatorname {deg } (#1)}
\def\pic#1.{\operatorname {Pic}\,(#1)}
\def\pico#1.{\operatorname{Pic}^0(#1)}
\def\picg#1.{\operatorname {Pic}^G(#1)}
\def\ner#1.{NS (#1)}
\def\rdown#1.{\llcorner#1\lrcorner}
\def\rup#1.{\ulcorner#1\urcorner}
\def\cone#1.{\operatorname {NE}(#1)}
\def\ccone#1.{\overline{\operatorname {NE}}(#1)}
\def\coef#1.{\frac{(#1-1)}{#1}}
\def\vit#1.{D_{\langle #1 \rangle}}
\def\mm#1.{\overline {M}_{0,#1}}
\def\H1#1.{H^1(#1,{\ring #1.})}
\def\ac#1.{\overline {\Bbb F}_{#1}}
\def\mgn#1.#2.{\overline {M}_{#1,#2}}
\def\ilist#1.{{#1}_1,{#1}_2,\dots}

\def\adj#1.{\frac {(#1-1)}{#1}}


\def\list#1.#2.{{#1}_1,{#1}_2,\dots,{#1}_{#2}}
\def\omitlist#1.#2.{{#1}_1,{#1}_2,\dots,\hat {{#1}_i}, \dots, {#1}_{#2}}
\def\omitlist0#1.#2.{{#1}_0,{#1}_1,\dots,\hat {{#1}_i}, \dots, {#1}_{#2}}
\def\loc#1.#2.{\Cal O_{#1,#2}}
\def\fderiv#1.#2.{\frac {\partial #1}{\partial #2}}
\def\map#1.#2.{#1 \longrightarrow #2}
\def\rmap#1.#2.{#1 \dasharrow #2}
\def\emb#1.#2.{#1 \hookrightarrow #2}
\def\non#1.#2.{\text {Spec }#1[\epsilon]/(\epsilon)^{#2}}
\def\Hi#1.#2.{\text {Hilb}^{#1}(#2)}
\def\sym#1.#2.{\operatorname {Sym}^{#1}(#2)}
\def\Hb#1.#2.{\text {Hilb}_{#1}(#2)}
\def\Hm#1.#2.{\Hom_{#1}(#2)}
\def\prd#1.#2.{{#1}_1\cdot {#1}_2\cdots {#1}_{#2}}
\def\Bl #1.#2.{\operatorname {Bl}_{#1}#2}
\def\pl #1.#2.{#1^{\otimes #2}}


\def\alist#1.#2.#3.{#1_1 #2 #1_2 #2\dots #2 #1_{#3}}
\def\zlist#1.#2.#3.{#1_0 #2 #1_1 #2\dots #2 #1_{#3}}
\def\lmap#1.#2.#3.{#1 \overset #2\to \longrightarrow #3}
\def\ses#1.#2.#3.{0\longrightarrow #1 \longrightarrow #2 \longrightarrow #3 
\longrightarrow 0}
\def\les#1.#2.#3.{0\longrightarrow #1 \longrightarrow #2 \longrightarrow #3}
\def\es#1.#2.#3.{#1 \longrightarrow #2 \longrightarrow #3}
\def\Hi#1.#2.#3.{\text {Hilb}^{#1}_{#2}(#3)}


\def\Hom{\operatorname{Hom}}

\def\deg{\operatorname{deg}}

\def\lcs{\operatorname{LLC}}


\def\e{\Cal E}

\def\e1{E_1}
\def\e2{E_2}

\def\ds{\displaystyle}

\def\qle{\sim_{\Bbb Q}}

\NoBlackBoxes
\topmatter
\title
Boundedness of log terminal Fano pairs of bounded index
\endtitle
  
\author
 James M\raise 1.6pt \hbox{\text {\smc c}}Kernan 
\endauthor  
\address
Department of Mathematics, University of California at Santa Barbara, Santa Barbara, CA 93106
\endaddress
\endtopmatter

\head \S 1 Introduction and statement of results\endhead

 A fundamental problem in classifying varieties is to determine natural subsets whose
moduli is bounded.  The difficulty of this problem is partially measured by the behaviour
of the canonical class of the variety.  Three extreme cases are especially of interest:
either the canonical class is ample, that is the variety is of general type, or the
canonical class is trivial, for example the variety is an abelian variety or Calabi-Yau,
or minus the canonical class is ample, that is the variety is Fano.

 Let us first suppose that the variety is smooth.  It is well known that varieties of
general type do not form a bounded family unless one fixes some invariants.  For example a
curve of genus $g$ is of general type iff $g\geq 2$.  Typically even varieties with
trivial canonical class do not form bounded families; for example the moduli space of
abelian varieties has infinitely many components corresponding to the type of
polarisation.  Even though the same is true of K3 surfaces, to the best of the author's
knowledge it is not known whether the family of Calabi-Yau threefolds is bounded.  The
picture for smooth Fano varieties however is much brighter; indeed there is a complete
classification up to dimension three.  Such a classification in general seem unfeasible,
but on the other hand it was proved by Koll\'ar, Miyaoka and Mori \cite{KMM92c} that smooth
Fano varieties of fixed dimension form a bounded family.

 One reason for focusing on the three extreme cases is that roughly speaking it is
expected that up to birational equivalence, any variety is either of general type or
admits a fibration to another variety, whose fibres are either Fano or have trivial
canonical class.  However one can only achieve this birational factorisation if one allows
singularities.  For example for a surface of general type one needs to contract all the
$-2$-curves on the surface to make the canonical divisor ample.

 Unfortunately it is no longer true that singular Fano varieties form a bounded family.
For example cones over a rational normal curve of degree $d$ form an unbounded family of
Fano surfaces; to form a bounded family one needs to impose some restrictions on the
singularities.  One natural restriction to impose is that some fixed multiple of the
canonical class is Cartier.  The smallest such positive integer we will call the index.
(Note that the index is also used to designate other invariants in this context;
unfortunately there does not seem to be a standard notation.)  However even this is not
enough, since cones over elliptic curves have index one, yet form an unbounded family of
Fano surfaces, simply varying the degree as before.

 One natural class of singularities to work with are kawamata log terminal singularities.
They form a large class of singularities which are closed under many standard geometric
operations; for example the finite quotient or unramified cover of any kawamata log
terminal singularity is kawamata log terminal.  Moreover they seem to be the largest class
where one might expect boundedness; for example cones over an elliptic curve are log
canonical but not kawamata log terminal.  Fortunately to achieve the factorisation alluded
to above it is expected that one only needs to allow the presence of canonical
singularities which are automatically kawamata log terminal.

 It is also natural to enlarge the category we are working in and include a boundary
divisor with fractional coefficients.  These appear naturally when one takes a quotient
and there is ramification in codimension one or when one passes to a resolution.  For
example suppose $S$ is the cone over a rational normal curve of degree $d$.  Then the
minimal resolution $T$ is isomorphic to $\Hz n.$, the unique $\pr 1.$-bundle over $\pr 1.$
with a section $C$ of self-intersection $-n$ and $-K_S-\adj n.C$ is ample.  The morphism
$\map T.S.$ simply contracts $C$.

\definition{1.1 Definition} We say that the pair $(X,\Delta)$ is a 
{\bf log pair} if $X$ is a normal variety and $\Delta$ is a $\Bbb Q$-divisor such that
$K_X+\Delta$ is $\Bbb Q$-Cartier.

 If $\pi\:\map Y.X.$ is a morphism the {\bf log pullback} $\Gamma$ of $\Delta$ is defined by the formula
$$
K_Y+\Gamma=\pi^*(K_X+\Delta).
$$

  We say that $\Delta$ is a {\bf sub-boundary} if all of its coefficients are less than one;  if 
in addition the coefficients of $\Delta$ are all positive then we say that $\Delta$ is a 
{\bf boundary}.  

 We say that the pair $(X,\Delta)$ is {\bf kawamata log terminal} if $\rdown \Delta.$ is
empty and the log pullback is a sub-boundary, for any birational morphism. 
\enddefinition

 It is impractical to give all the standard definitions and results of higher
dimensional geometry; the reader is referred to the excellent introduction of
\cite{KM98} for those definitions we have omitted.  We work over an algebraically 
closed field of characteristic zero.  We prove the following conjecture of Batryev:

\proclaim{1.2 Theorem} Let $r$ and $n$ be integers.  Then the family of all log 
pairs $(X,\Delta)$ such that
\roster 
\item $X$ has dimension $n$, 
\item $-(K_X+\Delta)$ is ample, 
\item $K_X+\Delta$ is kawamata log terminal, and 
\item $r(K_X+\Delta)$ is Cartier 
\endroster 
is bounded. 
\endproclaim

 As pointed out above, (1.2) includes Fano varieties with quotient singularities.
To prove (1.2), by a result of Koll\'ar \cite{Kollar93b}, it suffices to bound the
top self-intersection of $-(K_X+\Delta)$, which we will refer to as the degree $d$.
Thus (1.2) follows from

\proclaim{1.3 Theorem} Let $r$ and $n$ be integers.  Then there is a 
real number $M$ such that if we have a log pair $(X,\Delta)$ where
\roster 
\item $X$ has dimension $n$, 
\item $-(K_X+\Delta)$ is big and nef,  
\item $K_X+\Delta$ is kawamata log terminal, and 
\item $r(K_X+\Delta)$ is Cartier, 
\endroster 
then $d<M$.  
\endproclaim

 One interesting feature of smooth Fano varieties is that they have trivial fundamental
group.  The same is true of Fano varieties with kawamata log terminal singularities but
actually in this case a more natural and more general question is to consider the
fundamental group of the complement of the singular points.  In general this group is
larger but it is believed that the fundamental group ought still to be finite.  Using
(1.3) we are able to prove that there is a bound on the maximum degree of a
finite cover so that at least the algebraic fundamental group is finite.  

\proclaim{1.4 Corollary} Let $(X,\Delta)$ be a kawamata log terminal log 
pair such that $-(K_X+\Delta)$ is big and nef.  Let $U$ be any open subset of $X$ whose
complement is of codimension at least two.

 Then the algebraic fundamental group of $U$ is finite.  
\endproclaim

 Problems centered around boundedness of Fano varieties have received considerable
attention.  Nikulin \cite{Nikulin89} proved (1.2) in dimension two and Borisov
\cite{Borisov99} proved (1.2) in dimension three.  Nadel
\cite{Nadel91} proved the boundedness of smooth Fano varieties with Picard number one.  As
previously mentioned Koll\'ar, Miyaoka and Mori \cite{KMM92c} proved the boundedness of
smooth Fano varieties in every dimension.  Recently Clemens and Ran \cite{CR00} proved the
boundedness of a restricted class of Fano varieties with Picard number one.

 Most proofs begin with an idea that goes back to Fano.  If the degree is very large then
by Riemann Roch for any point $p$ belonging to the smooth locus we may find a $\Bbb
Q$-divisor $B$ $\Bbb Q$-linearly equivalent to $-K_X$ with very high multiplicity at $p$.
The question is how to use the existence of $B$ to obtain a contradiction.

 \cite{KMM92c} and \cite{Campana93} use reduction modulo $p$ to produce chains of rational
curves of small degree connecting any two points.  This easily bounds the multiplicity of 
$B$ and hence the degree.  

 Borisov \cite{Borisov99} uses a similar argument.  Unfortunately, due to the presence of
singularities, this is extremely delicate and it seems hard to generalise this method to
higher dimensions.

 Nadel \cite{Nadel91} uses the fact that through any point of $X$ there is a rational
curve of low degree.  If this curve intersects $B$ at $p$ then it must lie in $B$.  He
then produces a covering family of varieties $V_t$ with the same property, namely if a
rational curve of low degree meets $V_t$ then in fact it must lie in $V_t$.  It is easy to
see that this cannot happen when the Picard number is one.  The locus $V_t$ is essentially
the locus where the divisor $B$ has the same high multiplicity.  

 Clemens and Ran \cite{CR00} show that certain sheaves of differential operators are
semi-positive.  Note that the condition that $B$ has at least a given multiplicity at a
point is equivalent to saying that any derivative, of a local defining equation, of order
up to the multiplicity vanishes.  $B$ then corresponds to a quotient of negative degree of
some sheaf of differential operators of high degree on a general curve passing through $p$. 
This clearly contradicts semi-positivity.

  We use entirely different ideas.  The main ideas of our proof are derived from the
$X$-method, which was developed to prove the Cone Theorem.  We make use of some of the
sophisticated ideas introduced, especially by Kawamata, in an attempt to solve Fujita's
conjecture.  Another key feature of our proof is to use the connectedness of the locus of
log canonical singularities.  (This phenomena was first observed by Shokurov in the case
of surfaces.  Koll\'ar proved that connectedness holds in general as a slick application of
Kawamata-Viehweg vanishing.)  As such, just as with the proof of Clemens and Ran
\cite{CR00}, we do not use reduction modulo $p$.  In particular we give the first proof of
boundedness of smooth Fano varieties without using reduction modulo $p$.

 Instead of producing chains of rational curves of low degree, we produce covering chains
of positive dimensional Fano varieties $V_t$ of low degree.  Unfortunately we are unable
to smooth these varieties, as one can smooth chains of rational curves, see for example
\cite{Kollar96}, and these chains only exist if the degree is sufficiently large.
Moreover it is far more problematic to show directly that the relevant intersection
numbers are bounded, due to problems of excess intersection which do not arise when one
intersects a divisor and a curve.  On the other hand these subvarieties are log canonical
centres of appropriate log divisors and as such they have a very rich geometry.  More
precisely

\proclaim{1.5 Theorem} Let $(X,\Delta)$ be a log pair where $X$ is 
projective of dimension $n>1$, $\Delta$ is an effective divisor and $-(K_X+\Delta)$ is big
and nef.
\roster 
\item If $d> (n!)^n$ then there is a subvariety of the Hilbert scheme of $X$ such that if 
$f\:\map Y.B.$ is the normalisation of the universal family and $\pi\:\map Y.X.$ is the
natural morphism, then $\pi$ is birational, $f$ is a contraction morphism and the fibres
of $f$ are Fano varieties and so rationally connected.  
\item If further $d\geq (2^nn!)^n$ then the fibres of $f$ have degree at most $(2^{k+1}k!)^k$.  
\item If $r(K_X+\Delta)$ is Cartier for some positive integer $r$ and 
$$
d\geq (2n)^n\left (2^nr(n-1)!\right)^{n(n-1)}
$$
then the log pullback $\Gamma$ of $\Delta$ is effective in a neighbourhood of the generic
fibre of $f$.
\item If $K_X+\Delta$ is kawamata log terminal of log discrepancy at least $\epsilon$ and 
$$
d\geq \left(\frac {2n}{\epsilon}\right)^n\left (2^nr(n-1)!\right)^{n(n-1)}
$$
then $\pi$ is small in a neighbourhood of the generic fibre of $f$.
\endroster 
\endproclaim

 A key part of the proof of (1.5) is that we place almost no restrictions on
the type of singularities of the pair $(X,\Delta)$.  Furthermore, even if one were to
start with a smooth variety $X$, with empty divisor $\Delta$, then inductively one needs
to deal with the case of an arbitrary log pair, so working at this level of generality is
crucial to the proof.

   It is not so hard to prove (1.3) using (1.5).  For example it is
easy to use (1.5) to obtain an explicit bound on the degree when the Picard
number is one.

\proclaim{1.6 Corollary}  Let $(X,\Delta)$ be a kawamata log terminal 
log pair.  If $X$ is $\Bbb Q$-factorial of Picard number one and $-r(K_X+\Delta)$ is
Cartier and ample, then the degree is at most
$$
(2nr)^n\left (2^nr(n-1)!\right)^{n(n-1)}.
$$
\endproclaim

 I also hope that (1.2) will be the first (and very small) step towards an
eventual proof of a very interesting conjecture, due independently to Alexeev and Borisov:

\proclaim{1.7 Conjecture} Fix an integer $n$ and a positive real number $\epsilon$.  

 Then the family of all varieties $X$ such that 
\roster 
\item the log discrepancy of $X$ is at least $\epsilon$ and
\item $-K_X$ is ample,
\endroster 
is bounded. 
\endproclaim

 Indeed by (1.2) it suffices to bound the index of $K_X$.  In fact more is true;
the proof of (1.3) shows that if (1.7) is true in dimension $n-1$
then the degree of $-K_X$ in dimension $n$ is bounded in terms of the log discrepancy.

Alexeev \cite{Alexeev94} proved (1.7) in dimension two and Borisov and Borisov
\cite{BB92} proved (1.7) for toric varieties.  Kawamata \cite{Kawamata92a}
proved (1.7) for Fano threefolds with Picard number one and terminal
singularities and Koll\'ar, Miyaoka, Mori and Takagi \cite{KMMT00} proved that all Fano
threefolds with canonical singularities are bounded.  Despite this (1.3) has not even been
resolved in dimension three.  It is also interesting to note that recently Lin has proved
that one cannot drop the reference to the log discrepancy $\epsilon$ in (1.7)
even in the case where one replaces bounded by birationally bounded, see \cite{Lin01}.

Note that if a threefold $X$ is uniruled (that is $X$ is covered by rational curves) then,
as a consequence of the MMP in dimension three, $X$ is birational to a Mori fibre space,
$\pi\:\map Y.S.$.  Quite surprisingly it turns out that many Fano varieties are
birationally rigid in the sense that there is only one Mori fibre structure $\pi$, up to
birational equivalence.  In higher dimensions the existence of the MMP has not yet been
established.  We adopt the following ad hoc definition of birationally rigidity in the
general case.

\definition{1.8 Definition} We will say that a variety $X$ is {\bf birationally rigid} 
if for any two covering families of curves $C_i$, such that $-K_X\cdot C_i<0$, then 
$C_1$ is numerically equivalent to a multiple of $C_2$.  
\enddefinition

 Note that this definition does indeed imply the more usual one, at least in dimension
three, by a result of Batryev \cite{Batryev89} .  With this said, (1.5) has
the following unexpected consequence

\proclaim{1.9 Corollary} Let $(X,\Delta)$ be a log pair where $X$ is 
projective of dimension $n>1$, $\Delta$ is an effective divisor and $-(K_X+\Delta)$ is big
and nef.
 
 If $d>(n!)^n$ then $X$ is not birationally rigid.  
\endproclaim

One wonders if there is an even stronger result than (1.9)

\proclaim{1.10 Conjecture} Fix an integer $n$.  Then the family of 
log pairs $(X,\Delta)$ such that 
\roster 
\item $K_X+\Delta$ is kawamata log terminal, 
\item $-(K_X+\Delta)$ is ample, and 
\item $X$ is birationally rigid 
\endroster 
is bounded.   
\endproclaim

 (1.10) provides an obvious strategy to prove (1.7).  Note also that
(1.10) is only of interest when the Picard number is one.  It might also be
interesting to study covering families of tigers in other contexts especially when the
degree is not necessarily large.  Indeed there are many analogies between the geometry of
covering families of tigers and the geometry of covering families of rational curves, some
of which we draw out below.  Finally it might be of interest to try to connect questions
of boundedness with the phenomena of Mirror symmetry.  Typically Mirror symmetry
associates to one moduli space of Calabi-Yaus another such space.  However in some
degenerate cases it switches Calabi-Yaus and Fano varieties.  Moreover many of the
constructions of Mirror symmetry implicitly involve Fano varieties.   

{\bf Sketch of Proof of (1.3):} To navigate through the proof (1.3), let
us use the example of cones over a rational normal curve of degree $d$.  The first step is
to pick up a covering family of tigers, see \S 3.  As pointed out above this corresponds
to finding a divisor $\Delta_p$ passing through $p$ of very large multiplicity.  The whole
point of covering families of tigers, is that we only really need to consider what happens
at a general point $p$ of $X$ (see for example (3.2)) and so in practice when we
need to apply adjunction or inversion of adjunction, (5.2) and
(6.1) we can treat the pair $(X,\Delta)$ as though it were log smooth.
Roughly speaking the locus $V_t$ is where $\Delta_t$ has very large multiplicity.  In this
sense a covering family of tigers $V_t$ is somewhat similar to the covering family
constructed by Nadel.

 In the case of a cone there are three obvious candidates,
\roster 
\item The family $V_t$ of lines, where $\Delta_t=V_t$.  
\item A pencil of hyperplane sections $W_s$, where $\Gamma_s=W_s$.
\item The family of points $p_r$ on $S$, where $\Theta_r=\Delta_t+\Gamma_s$.    
\endroster 

 As $-K_S=(d+2)\Delta_t$, the weight of the first family is $d+2$ and as $-K_S=\frac
{d+2}d\Gamma_s$ the weight of the second family is $\frac {d+2}d$.  The third family has
weight $\frac {d+2}{d+1}$, the harmonic sum of the first two.
 
The dimension of the first two families is obviously one and of the last zero.  In
practice then we never choose the second or third family as their weights are not large
enough.

 The next step is to refine the covering family $V_t$ so that it forms a birational
family, see \S 4.  In the case of a cone this is automatic.  In the general case it means
passing to a covering family of maximal dimension.  The idea is to average over all
divisors $\Delta_t$ passing through $p$, see (4.2).  This is somewhat
analogous to growing chains of rational curves.  The eats up some of the weight and we end
up with a birational family of much smaller weight.

 The third step is to decrease the degree of $V_t$, see \S 5.  Again in the case of a cone
this is automatic.  In the general case it means passing to a birational family of minimal
dimension.  The idea is that if the degree of $V_t$ is very large then we can lift a
birational family from $V_t$ to the whole variety, see (5.3).  This is quite
similar to bend and break for curves.  This eats up more of the weight and we end up with
a birational family of smaller weight still.

 Comparing the first family with the second it becomes clearer why we prefer the first;
the degree of $V_t$ is bounded whilst the degree of $W_t$ is arbitrarily large.

 A key observation is that in the Fano case the family must be non-trivial (that is $V_t$
is not of dimension zero) provided the weight is greater than two, see (3.4); this
is a simple consequence of connectedness.  This fact and the fact that we have to refine
our covering family twice is the reason for the rather large numbers that appear in
(1.5).

 In the final step we want to connect the log discrepancy to the degree of the fibres, see
\S 6.  To do this consider the natural diagram determined by a covering family of tigers
(3.1.1).  In the case of a cone we get
$$
\diagram
V_t  & \subset & T           &   \rTo^\pi & S,  \\
\dTo &         & \dTo^f      &            &     \\
t    & \in     & \pr 1.      &            &
\enddiagram
$$
where $T$ is isomorphic to $\Hz d.$.  The fibres of $f$ are isomorphic to $V_t$.  Thus the
fibres of $f$ are of bounded degree.  On the other hand if we pick two fibres of $f$,
$V_s$ and $V_t$ then they are not connected but their images in $X$ are connected.  Again
by connectedness it follows that there is a log canonical centre which connects them.  In
our case this is the exceptional divisor $E$ of $\pi$.  As the log discrepancy of $S$ is
$2/d$ it follows that the coefficient of $E$ with respect to $\Delta_t+\Delta_s$ is at
least $2/d$, that is, by symmetry, the coefficient of $E$ with respect to $\Delta_t$ is at
least $1/d$.  (In fact it is exactly $1/d$ corresponding to the fact that
$K_S+\Delta_s+\Delta_t$ is log canonical.)  The general case is not much harder, see
(6.2).

 By definition of the weight, we can find a $\Bbb Q$-divisor $B\qle -(K_X+\Delta)$ such
that the coefficient of $E$ with respect to $B$ is at least $w/d$.  Now $E$ is a section
of $f$ and $\pi^*(K_X+B)$ restricts to a trivial divisor on the general fibre of $f$.
Thus $w/d\leq 2$ and so $w\leq 2d$.  In turn this bounds the degree of $-(K_X+\Delta)$.
This gives a bound on the degree of $S$ only in terms of the log discrepancy.  In general
if $\pi$ is not small, a similar proof applies, see (6.4).  The general case is
similar but more subtle as in general the log canonical centre connecting $V_s$ and $V_t$
in $Y$ need not be a divisor, see (6.5).  This is the only place where the
argument is not effective.  

{\bf Thanks}: I would like to thank Florin Ambro for many useful conversations.  I would
also like to thank Professor Shokurov for some helpful comments.  The majority of this
work was completed whilst visiting Tokyo University.  I would like to thank Professor
Kawamata for his extremely kind hospitality during my stay in Tokyo.  This paper uses the commutative diagram package written in \TeX\ by Paul Taylor.

\head \S 2 Calculus of log canonical centres\endhead

 We refer the reader to \cite{Kollaretal} for the various notions of log terminal and log
canonical.  We first recall some useful definitions due to Kawamata \cite{Kawamata97} and
collect together some known results that we require for the proof of (1.2) and
(1.5).  

\definition{2.1 Definition}  Let $(X,\Delta)$ be a log pair.  A subvariety $V$ 
of $X$ is called a {\bf log canonical centre} if it is the image of a divisor of log
discrepancy at most zero.  A {\bf log canonical place} is a valuation corresponding to a divisor
of log discrepancy at most zero.  The {\bf log canonical locus} $\lcs(X,\Delta)$ of the
pair $(X,\Delta)$ is the union of the log canonical centers.  We will say that a log
canonical centre is {\bf pure} if $K_X+\Delta$ is log canonical at the generic point of
$V$.  If in addition there is a unique log canonical place lying over the generic point 
of $V$ we will say that $V$ is {\bf exceptional}.  
\enddefinition

 Log canonical centres enjoy many special properties.  We collect together six of the 
most important for us.  

\proclaim{2.2 Lemma} (Calculus of Log canonical centres) Let $(X,\Delta)$ 
be a log pair.    
\roster 
\item (Connectedness) If $f\:\map X.Z.$ is a morphism with connected fibres such that the image of 
every component of $\Delta$ with negative coefficient is of codimension at least two in $Z$ and 
$-(K_X+\Delta)$ is big and nef over $Z$ then the intersection of $\lcs(X,\Delta)$ with 
every fibre is connected.  
\item  $\lcs(X,\Delta)$ is a finite union of log canonical centres.
\item (Convexity) If $V$ is a log canonical centre of $\ds K_X+\Delta+\sum_{i=1}^d \Delta_i$, 
then there is an $i$ such that $V$ is a log canonical centre of $K_X+\Delta+d\Delta_i$.
\endroster
 Now suppose that in addition $\Delta$ is effective and that there is a boundary
$\Delta'\leq \Delta$ such that $K_X+\Delta'$ is kawamata log terminal (this always holds
for example if $X$ is log terminal and $\Bbb Q$-factorial).
\roster
\item[4] The intersection of two log canonical centres is a union of log canonical centres.  
\item If $p$ is a point of $\lcs(X,\Delta)$ then there is a unique minimal log canonical centre
$V$ containing $p$.
\item Suppose $V$ is a pure log canonical centre of $(X,\Delta)$ and $D$ is a $\Bbb Q$-Cartier 
divisor of the form $A+E$, where $A$ is ample and $E$ is effective.  If $V$ is not
contained in the support of $E$, then for all $\epsilon$ sufficiently small, there exists a
$\Bbb Q$-divisor $\Gamma\qle aD$, such that $V$ is an extremal log canonical centre of
$(X,\Delta-\epsilon(\Delta-\Delta')+\Gamma)$, where $a$ approaches zero as $\epsilon$
approaches zero.
\endroster 
\endproclaim
\demo{Proof} (1) follows immediately from (17.4) of \cite{Kollaretal}.  (2) follows by taking
a resolution of the pair $(X,\Delta)$ and observing that $\lcs(X,\Delta)$ is the image of
the exceptional divisors of log discrepancy less than or equal to zero.  Note that the
convex linear combination of divisors which are kawamata log terminal at the generic point
of $V$ is kawamata log terminal at the generic point of $V$.  As $\sum \Delta_i$ is a
convex linear combination of the divisors $d\Delta_i$, (3) follows.  (4) is (1.5) of
\cite{Kawamata97} and (5) is immediate from (4).  (6) is observed in \cite{Kawamata98}. \qed\enddemo

 In a recent preprint of Ambro \cite{Ambro01} it is proved that (2.2.4) and
(2.2.5) hold, even without the unnatural hypothesis on the existence of 
$\Delta'$.

\head \S 3 Covering families of tigers\endhead

\definition{3.1 Definition} Let $(X,\Delta)$ be a log pair, where $X$ is 
projective and let $D$ be a $\Bbb Q$-Cartier divisor. We say that pairs of the form
$(\Delta_t,V_t)$ form a {\bf covering family of tigers of dimension $k$ and weight $w$} if
\roster 
\item There is a projective morphism $f\:\map Y.B.$ of normal projective varieties and an open subset 
$U$ of $B$ such that the fibre of $\pi$ over $t\in U$ is $V_t$.  
\item There is a morphism of $B$ to the Hilbert scheme of $X$ such that $B$ is the normalisation of its
image and $f$ is obtained by taking the normalisation of the universal family.  
\item If $\pi\:\map Y.X.$ is the natural morphism then $\pi(V_t)$ is a pure log canonical centre of 
$K_X+\Delta+\Delta_t$. 
\item $\pi$ is dominant. 
\item $\pi$ is finite.  
\item $D\qle w\Delta_t$, where $\Delta_t$ is effective. 
\item the dimension of $V_t$ is $k$.   
\endroster

 We will say that the family is a {\bf birational family of tigers} if 
\roster
\item[8] $\pi\:\map Y.X.$ is birational. 
\endroster 

 Further if 
\roster
\item[9]
$$
w> 2^n(n-k)!
$$
\endroster 
then we will say that $(\Delta_t,V_t)$ forms a {\bf covering family of tigers of large
weight}.

 Finally if 
\roster
\item[10] 
$$
w> 2^{k+1}
$$
\endroster
and the family is birational then we will say that $(\Delta_t,V_t)$ forms a {\bf
birational family of tigers of large weight}.
\enddefinition

 The following commutative diagram is perhaps worth a couple of words:
$$
\diagram
V_t  & \subset & Y      &   \rTo^\pi & X.  \\
\dTo &         & \dTo^f &            &     \\
t    & \in     & B      &            &
\enddiagram
\tag 3.1.1
$$

 Note that $\pi|_{V_t}\:\map V_t.\pi(V_t).$ is finite and birational.  We will often abuse
notation and identify $V_t$ with its image.  Clearly (3.1.6) may be rewritten
as $\Delta_t\qle \frac 1wD$.  In particular the weight of a sum of divisors is the
harmonic sum (that is the reciprocal of the sum of the reciprocals) of the weights.  The
whole point of (3.1) is that covering family of tigers enjoy some very
special properties.  When we apply (3.1), we will always take
$D=-(K_X+\Delta)$.  As the general case is no harder and with a view towards other
possible applications, we take $D$ to be arbitrary.

 We first show that covering families of tigers exist under weak conditions and at the
same time justify the notation of (3.1).  
 
 Given a topological space $X$, we will say that a subset $P$ is countably dense if $P$ is
not contained in the union of countably many closed subsets of $X$.  Note that if we decompose
as $P$ a countable union 
$$
P=\bigcup_{m\in\Bbb N} P_m
$$
then at least one $P_m$ is countably dense.  

\proclaim{3.2 Lemma} Let $(X,\Delta)$ be a log pair, where $X$ is projective
and let $D$ be a big $\Bbb Q$-Cartier divisor.  Let $w$ a positive rational number.  Let
$P$ be a countably dense subset of $X$.  Suppose that for every point $p\in P$ we may find
a pair $(\Delta_p,V_p)$ such that $p\in V_p$ and $V_p$ is a pure log canonical centre of
$K_X+\Delta+\Delta_p$ where $\Delta_p\qle D/w_p$ for some $w_p>w$.

 Then we may find a covering family of tigers of weight $w$, together with a countably
dense subset $Q$ of $P$ such that for all $p\in Q$, $V_p$ is a fibre of $\pi$.
\endproclaim
\demo{Proof} As $D$ is big, it is of the form $A+E$ where $A$ is ample and $E$ is effective.  
By (2.2.6), possibly passing to a subset of $P$, we may assume that $V_p$ is
an exceptional log canonical centre.  Adding on a multiple of $D$ we may assume that
$w_p=w$.

 We may find an integer $m$ so that $mD/w$ is integral.  Moreover for each $p\in P$ there
is an integer $m_p$ such that $m_p\Delta$ is linearly equivalent to $m_pD/w$.  Hence using
the observation above, we may assume there is a countably dense subset $Q$ of $P$ such
that for all $p\in Q$, $m\Delta_p\in |mD/w|$.  Let $B$ be the closure inside this linear
system of the points corresponding to $m\Delta_p$ for $p$ in $Q$.  Replace $B$ by an
irreducible component that contains a countably dense set of points $Q$ in $X$ and let
$f\:\map H.B.$ be the universal family.  Pick a log resolution of the generic fibre of the
universal family and extend this to an embedded resolution over an open subset $U$ of $B$.
By assumption there is a unique exceptional divisor of log discrepancy zero over the
generic point corresponding to the log canonical centres $V_p$, for $p\in Q$.

 Thus possibly taking a finite cover of $B$ and passing to an open subset of $U$ we may
assume that there is a morphism $f\:\map Y.B.$ whose fibre $V_t$ over $t\in U$ is a log
canonical centre of $K_X+\Delta+\Delta_t$.  Possibly passing to an open subset of $U$ we
may assume that $f$ is flat, so that $U$ maps to the Hilbert scheme.  Replacing $B$ by the
normalisation of the closure of the image of $U$ in the Hilbert scheme and $Y$ by the
normalisation of the pullback of the universal family we may assume that
(3.1.1-2) hold.  Cutting by hyperplanes in $B$ we may assume that $f$ is
finite.  The result is now clear. \qed\enddemo

We will use (3.2) repeatedly, often without comment.

\proclaim{3.3 Lemma} Let $(X,\Delta)$ be a log pair, where $X$ is a 
projective variety of dimension $n$ and suppose that $D$ is big and nef.  Let $w$ be a
positive rational number and let $p$ and $q$ be two points of $X$ contained in the smooth
locus of $X$ but not contained in a component of $\Delta$ of negative coefficient.  Let
$d$ be the degree of $D$.
\roster 
\item If $d>(wn)^n$ then we may find an effective divisor $\Gamma\qle D/w $ such that 
$\lcs(X,\Delta+\Gamma)$ contains $p$.
\item If $d>(wn)^n$ then we may find a covering family of tigers of weight $w$. 
\item If $d>(2^nn!)^n$ then we may find a covering family of tigers of large weight.
\endroster
 Now assume further that $\Delta$ is effective and that $-(K_X+\Delta+D)$ is big and nef.  
\roster
\item[4] If $d>(2n)^n$ then there is an effective divisor $\Gamma\qle D$ and 
a chain $\list V.k.$ of log canonical centres of $K_X+\Delta+\Gamma$ such that $p\in V_1$,
$q\in V_k$ and $V_i\cap V_{i+1}\neq\emptyset$, for $i\leq k-1$.
\endroster 
\endproclaim
\demo{Proof} (1) is well known, see for example (7.1) of \cite{KM99} for a proof. 
(2) follows from (1) and (3.2).  (3) follows from (2).

 By (1) we may find a divisor $\Gamma_p\qle D/2$ (respectively $\Gamma_q$) such that
$K_X+\Delta+\Gamma_p$ is not kawamata log terminal at $p$ (respectively $q$).  Let
$\Gamma=\Gamma_p+\Gamma_q\qle D$.  Then $K_X+\Delta+\Gamma$ is not kawamata log terminal
at $p$ and $q$ and $-(K_X+\Delta+\Gamma)$ is big and nef.  (4) now follows from
(2.2.1).  \qed\enddemo

\proclaim{3.4 Lemma} Let $(X,\Delta)$ be a log pair, where $X$ is projective and 
$\Delta$ is effective and let $D$ be a $\Bbb Q$-Cartier divisor such that
$-(K_X+\Delta+D)$ is big and nef.  Let $(\Delta_t,V_t)$ be a covering family of tigers of
weight $w$ greater than two.

 Then $k>0$.  
\endproclaim
\demo{Proof} Suppose not.  Then for general $p_1$ and $p_2$ we may find $\Delta_1$ and
$\Delta_2$ such that $\Delta_i\qle D/w$ where $p_i$ is an exceptional log canonical centre
of $K_X+\Delta+\Delta_i$ and $w>2$.  As $p_1$ and $p_2$ are general it follows that
$\Delta_2$ does not contain $p_1$.  Thus $p_1$ is an exceptional log canonical centre of
$K_X+\Delta+\Delta_1+\Delta_2$ and hence $\lcs(X,\Delta+\Delta_1+\Delta_2)$ is not
connected.  On the other hand as $w>2$, $-(K_X+\Delta+\Delta_1+\Delta_2)$ is big and nef.
This contradicts (2.2.1) (or indeed (3.3.4)).  \qed\enddemo

\head \S 4 Covering families of maximal dimension\endhead

We need an easy Lemma that bounds, both from below and above, the multiplicity of a log
canonical divisor at a smooth point of a variety.  Even though we are only interested in
the case of effective divisors, for the purposes of induction it is expedient to allow
negative coefficients.

\proclaim{4.1 Lemma} Let $X$ be a smooth variety of dimension $n$.
Let $D$ be a normal crossings divisor, let $\Delta$ be a boundary and $\Gamma$ an
effective divisor.  Let $Z$ be a pure log canonical centre of $K_X+D+\Delta-\Gamma$ of
dimension $k$.  Let $\mu$ be the multiplicity of $\Delta$ at the generic point of $Z$
and let $a$ be the maximum coefficient of $D$.  Assume that $a<1$.  
\roster 
\item If $\Gamma$ is empty then $\mu\leq n-k$.  
\item $\mu>1-a$. 
\endroster 
\endproclaim
\demo{Proof} The result is local about the generic point of $Z$, and so we may as well 
assume that $X$ is affine.  Cutting by hyperplanes in $X$ we may assume that $Z$ is a
point $p$.  Let $\pi\:\map Y.X.$ denote the blow up of $X$ at $p$ with exceptional divisor
$E$.  Then the log discrepancy of $E$ with respect to $K_X+\Delta$ is $n-\mu$, which is at least
zero as $K_X+\Delta$ is log canonical.  Hence (1).

 On the other hand, by assumption there is an algebraic valuation $\nu$ of log discrepancy
zero with respect to $K_X+\Delta-\Gamma$.  $\nu$ determines a series of blows ups, each
with centre $\nu$, such that eventually the centre of $\nu$ is a divisor, see for example
(2.45) of \cite{KM98}.  Suppose that $\mu\leq 1-a$.  We have
$$
K_Y+D'+E=\pi^*(K_X+D)+bE,
$$
where $D'$ is the strict transform of $D$ and $b\geq (n-na)=n(1-a)$.  Thus the log
discrepancy of $E$ with respect to $K_X+D+\Delta-\Gamma$ is at least 
$$
n(1-a)-\mu\geq (n-1)(1-a)\geq (1-a).
$$  
Moreover the coefficient of $E$ is then at most $a$, $D'\cup E$ still has normal crossings
in a neighbourhood of the generic point of the centre of $\nu$, and the multiplicity of
the strict transform of $\Delta$ is at most the multiplicity of $\Delta$.  By induction on
the number of blow ups, it follows that $\nu$ has log discrepancy at least $1-a$, a
contradiction. \qed\enddemo

\proclaim{4.2 Lemma} Let $(X,\Delta)$ be a log pair, where $X$ is a normal 
projective variety and let $D$ be a big and nef divisor.  Let $(\Delta_t,V_t)$ be a
covering family of tigers of weight less than $w$ and dimension $k$.

 If $(\Delta_t,V_t)$ is not birational then we may find a covering family of tigers
$(\Gamma_s,W_s)$ of weight $w/(n-k)$ and dimension $l$, where either 
\roster 
\item $l>k$, or 
\item $l<k$ and $(\Gamma_s,W_s)$ is a birational family.  
\endroster 
\endproclaim
\demo{Proof} Let $d$ be the degree of $\pi$.  By assumption $d>1$.  Pick an open subset 
$U\subset X$ contained in the smooth locus of $X$ and the complement of the support of
$\Delta$ such that for all $p\in U$ there exist $\list t.d.$ points of $B$ such that $p\in
V_t$ iff $t=t_i$, $i=1\dots d$.  Set $\Delta_i=\Delta_{t_i}$ and $V_i=V_{t_i}$.  Let
$G=\sum \Delta_i$ and pick $\phi$ such that $K_X+\Delta+\phi G$ is maximally log canonical
at $p$.  As $p$ is a smooth point of $X$, by (4.1) we have
$$
n-k\geq \mu_p(\phi G)=\phi\sum _i\mu_p(\Delta_i)\geq d\phi.
$$

Let $W$ be the unique minimal log canonical centre of $K_X+\Delta+\phi G$ containing $p$.
There is an open subset of $\alist V.\cap.d.$ such that as we vary $p$ in this set, $W$
does not change.  Thus we may assume $\alist V.\cap.d.\subseteq W$.  If we have equality
then by (3.2) there is a small perturbation $\Gamma$ of $\phi G$ such that
$(\Gamma_s,W_s)$ forms a covering family of weight $w/(n-k)$, which is obviously
birational.

 Thus we may assume that $W\subsetneq V_1$.  By symmetry and (2.2.3) $W$ is a
log canonical centre of $K_X+d\phi \Delta_1$.  Pick $\lambda\leq d\phi$ so that $W$ is a
pure log canonical centre of $K_X+\Delta+\lambda\Delta_1$.  Now varying $p\in V_1$ it
follows that there is an irreducible component of $\lcs(X,\Delta+\lambda\Delta_1)$ that
contains $V_1$ as a proper subset.  By (2.2.2), we may assume that $V_1$ is a
proper subset of a log canonical centre $W_1$ of $\lcs(X,\Delta+\lambda\Delta_1)$.  It
follows by (3.2) that there is a small perturbation $\Gamma$ of
$\lambda\Delta_1$, so that $(\Gamma_s,W_s)$ is a covering family of tigers of weight
$w/(n-k)$ and greater dimension.  \qed\enddemo

\proclaim{4.3 Lemma} Let $(X,\Delta)$ be a log pair, where $X$ is a normal 
projective variety.  Let $D$ be a big and nef divisor.  Let $(\Delta_t,V_t)$ be a covering
family of tigers of large weight $w$ and maximal dimension $k$.
 
Then either 
\roster
\item $(\Delta_t,V_t)$ is a birational family of tigers of large weight, or   
\item there is a birational family $(\Gamma_s,W_s)$ of tigers of large weight and 
smaller dimension.  
\endroster
\endproclaim
\demo{Proof} By assumption $w>2^n(n-k)!\geq 2^{k+1}$.  If $(\Delta_t,V_t)$ is a birational 
family of tigers (1) holds and there is nothing to prove.  Otherwise by (4.2)
there is a covering family of tigers of weight $w'>2^n(n-k-1)!$.  The dimension of this
family cannot be larger than $k$, otherwise $(\Gamma_s,W_s)$ would be a covering family
$(\Gamma_s,W_s)$ of tigers of large weight and greater dimension, which contradicts the
maximality of $(\Delta_t,V_t)$.  By (4.2) the only other possibility is that
$(\Gamma_s,W_s)$ is a birational family of tigers of large weight and smaller dimension
and we have (2). \qed\enddemo

\head \S 5 Birational Families of minimal dimension\endhead

We need a form of adjunction and inversion of adjunction.  Conjecturally there ought to be
a strong form of adjunction for the normalisation of any log canonical centre.
Unfortunately inversion of adjunction is only known for divisors and adjunction is only
known for some special centres (for example minimal), see \cite{Kawamata98} and
\cite{Ambro01}.  Fortunately we only need adjunction and inversion of adjunction in very
special cases.  We state one form of adjunction to do with log canonical centres here.
The other form concerns birational families, see (6.1).

\proclaim{5.1 Lemma} Let $\pi\:\map Y.X.$ be a smooth morphism of 
smooth varieties.  Let $\Delta_1$ be a $\Bbb Q$-divisor on $X$ and let $\Gamma_1$ be the
pullback of $\Delta_1$ to $Y$.  Let $\Gamma_2$ be a boundary on $Y$, such that the 
support $B$ of $\Gamma_2$ dominates $X$ and $\pi|_B$ is smooth.  

 Then $(X,\Delta)$ is log canonical (respectively kawamata log terminal, etc.) iff
$(Y,\Gamma=\Gamma_1+\Gamma_2)$ is log canonical (respectively kawamata log terminal,
etc.).
\endproclaim
\demo{Proof} The property of being log canonical is local in the analytic topology. 
On the other hand, locally in the analytic topology, $Y$ is a product $X_1\times X_2$,
where $X_1$ is isomorphic to $X$ and $\Gamma_2$ is the pullback of a divisor $\Delta_2$
whose support is smooth, so that $\Gamma=\Gamma_1\times X_2+X_1\times \Gamma_2$.  But the
result is easy in this case, see for example (8.21) of \cite{Kollar95}.\qed\enddemo

\proclaim{5.2 Lemma}  Let $(X,\Delta)$ be a log pair where $X$ 
is projective and $\Delta$ is effective.  Suppose that $V$ is an exceptional log canonical
centre of $K_X+\Delta$.  Then there is an open subset $U$ of the smooth locus of $V$ with
the following property:

 For all divisors $\Theta$ on $X$, which do not contain the generic point of $V$ and
subvarieties $W$ of $V$ such that $W|_U$ is a pure log canonical centre of $K_U+\Theta|_U$
then $W$ is a pure log canonical centre of $K_X+\Delta+\Theta$.
\endproclaim
\demo{Proof} This result is local about the generic point of $V$ so we are free to replace
$X$ by any open set that contains the generic point of $V$.  The idea is to reduce to the
case of a divisor, when the result becomes an easy consequence of inversion of adjunction.
Pick a log resolution $\pi\:\map Y.X.$ of the pair $(X,\Delta)$ and let $\Gamma$ be the
log pullback of $\Delta$.  By assumption there is a unique divisor $E$ lying over $V$ of
log discrepancy zero.  Let $f\:\map E.V.$ be the restriction of $\pi$ to $E$.  Replacing
$X$ by an open subset, we may assume that $f$ and $V$ are both smooth, and that $K_V$ and
$K_X+\Delta$ are $\Bbb Q$-linearly equivalent to zero.  By adjunction we may write
$$
(K_Y+\Gamma)|_E=K_E+B,
$$
for some effective divisor $B$, where both sides are $\Bbb Q$-linearly equivalent to zero.
Possibly passing to a smaller open set we may assume that every component of $B$ dominates
$V$ and that $f|_B$ is smooth.  Then
$$
K_E+B=(K_Y+\Gamma)|_E=\pi^*(K_X+\Delta)|_E=f^*((K_X+\Delta)|_V)\qle f^*(K_V).
$$
Suppose that $W$ is a pure log canonical centre of $K_V+\Theta|_V$.  Set
$\Theta'=\pi^*\Theta$.  As $\Theta$ does not contain the generic point of $V$, $E$ is not
a component of $\Theta'$, so that $f^*(\Theta|_V)=\Theta'|_E$.  It follows by
(5.1) that $f^{-1}(W)$ is a pure log canonical centre of $K_E+\Theta'|_E$.  The
result now follows as in the proof of (17.1.1), (17.6) and (17.7) of \cite{Kollaretal}.
\qed\enddemo

\proclaim{5.3 Lemma} Let $(X,\Delta)$ be a log pair and let $D$ be 
a divisor of the form $A+E$ where $A$ is ample and $E$ is effective.  Let $(\Delta_t,V_t)$
be a covering family of tigers of weight greater than $w$ and dimension $k$.  Let $A_t$ be
the restriction of $A$ to $V_t$.  Suppose that there is an open set of $B$ such that for
all $t\in B$ we may find a covering family of tigers $(\Gamma_{t,s},W_{t,s})$ on $V_t$ of
weight, with respect to $A_t$, greater than $w'$.

 Then we may find a covering family of tigers $(\Gamma_s,W_s)$ of dimension less than $k$ 
and weight 
$$
w''=\frac 1w+\frac 1{w'}=\frac {ww'}{w+w'}.
$$
Further if both $(\Delta_t,V_t)$ and $(\Gamma_{t,s},W_{t,s})$ are birational families then 
so is $(\Gamma_s,W_s)$.  
\endproclaim
\demo{Proof} Pick $r$ so that $rA$ is Cartier and let $L=\ring X.(rA)$ be the corresponding 
line bundle.  Note that by Serre vanishing $H^1(X,\Cal I_V\otimes \pl L.m.)=0$ for $m$
large enough.  Hence we may lift $\Gamma_{t,s}$ to a $\Bbb Q$-divisor $\Theta_{t,s}$ on
$X$ $\Bbb Q$-linearly equivalent to the same multiple of $A$.  Adding on a multiple of $E$
we may assume that $\Theta_{t,s}$ is in fact a multiple of $D$.  Thus by
(5.2) and (3.2) applied to $(\Delta_t+\Theta_{t,s},W_{t,s})$ we
may find a covering family of tigers $(\Gamma_s,W_s)$ of weight $w''$.  The rest is clear.
\qed\enddemo

\proclaim{5.4 Lemma} Let $(X,\Delta)$ be a log pair and let $D$ be an 
ample divisor.  Let $(\Delta_t,V_t)$ be a birational family of tigers of large weight with
minimal dimension $k$.

 Then $D^k\cdot V_t$ is at most $(2^{k+1}k!)^k$. 
\endproclaim
\demo{Proof} Suppose that $D^k\cdot V_t>(2^{k+1}k!)^k$.  Since $D$ is big we may write $D\qle A+E$,
where $A$ is ample and $E$ is effective.  By Kleiman's criteria for ampleness, as $D$ is
nef, $D-\epsilon E$ is ample for any $0<\epsilon\leq 1$.  Thus we may assume that
$A^k\cdot V_t>(2^{k+1}k!)^k$.  Set $A_t=A|_{V_t}$.  Then $(1/2A_t)^k>(2^kk!)^k$.  By
(4.2) there is a birational family of tigers $(\Gamma_{t,s},W_{t,s})$ of
weight greater than $2^k$ for $K_{V_t}+\Gamma_t$ and $1/2A_t$, that is of weight greater
than $2^{k+1}$ for $A_t$.  By (5.3) it follows that there is a birational family
of tigers $(\Gamma_s,W_s)$ of dimension less than $k$ and weight greater than $2^k$, that
is of large weight, which contradicts the minimality of $(\Delta_t,V_t)$. \qed\enddemo

\head \S 6 A bound for the degree of a birational family in terms of the log discrepancy\endhead

 Here is the other version of adjunction we shall need, mentioned in \S 5.

\proclaim{6.1 Lemma} Let $X$ be a projective variety and suppose we are 
given a subvariety of the Hilbert scheme of $X$ such that if $f\:\map Y.B.$ is the
normalisation of the universal family and $\pi\:\map Y.X.$ is the natural morphism, then
$\pi$ is birational.  Let $\Delta$ be a $\Bbb Q$-divisor such that $K_X+\Delta$ is $\Bbb
Q$-Cartier.

 Then there is an open subset $U\subset B$ such that for all $t\in U$, the fibre $V_t$
of $f$ over $t$ is smooth in codimension one and we may find $\Gamma_t$ such that 
$$
(K_X+\Delta)|_{V_t}=K_{V_t}+\Gamma_t.
$$
 
 Moreover the pair $(V_t,\Gamma_t)$ has the following properties. 
\roster
\item $\Gamma_t$ is effective iff the log pullback $\Gamma$ of $\Delta$ is effective in a 
neighbourhood of the generic fibre of $f$.
\item If $\Gamma$ is effective in a neighbourhood of the generic fibre of $f$ and $K_X+\Delta$ 
is kawamata log terminal then $V_t$ is normal.
\item The log discrepancy of $(V_t,\Gamma_t)$ is at least the log discrepancy of $(X,\Delta)$.  
\endroster 
\endproclaim
\demo{Proof} It suffices to define $\Gamma_t$ and to prove (1-3) on $Y$, in a neighbourhood 
of the generic fibre of $f$.

The generic fibre of $f$ is certainly smooth in codimension one and so the fibres of $f$
are certainly smooth in codimension one over an open subset.  Thus over the same open
subset $K_X|_{V_t}=K_{V_t}$ and so we can define $\Gamma_t$ to simply be the restriction
of $\Gamma$ to $V_t$.  (1) is then clear.

 If $\Gamma_t$ and $K_X+\Delta$ is kawamata log terminal it follows that $K_Y+\Gamma$ is
kawamata log terminal of log discrepancy at least the log discrepancy of $K_X+\Delta$ in a
neighbourhood of the generic fibre of $f$.  In particular it follows that $Y$ is Cohen
Macaulay in a neighbourhood of the generic fibre of $f$ and so the fibres of $f$ are Cohen
Macaulay over an open subset.  In particular $V_t$ is normal.  Hence (2).

  Picking a resolution of the generic fibre and possibly passing to a smaller open set,
(3) follows easily. \qed\enddemo

\proclaim{6.2 Lemma} Let $(X,\Delta)$ be a log pair, where $X$ is projective.  
Let $D$ be any $\Bbb Q$-Cartier divisor and $(\Delta_t,V_t)$ a birational family of
tigers.  Set $\Gamma_t=\pi^*\Delta_t$.
  
 Then there is an open subset $U$ of $B$ such that if $t\in U$ then every component $Z$ of
the $\pi$-exceptional locus which dominates $B$ is a log canonical centre of
$K_X+\Gamma+2\Gamma_t$.  In particular every $\pi$-exceptional divisor has log discrepancy
at most zero with respect to $K_X+\Gamma+2\Gamma_t$.
\endproclaim
\demo{Proof} Pick $s\neq t\in B$.  Then $V_s$ and $V_t$ are log canonical centres of 
$K_Y+\Gamma+\Gamma_s+\Gamma_t$.  Thus (2.2.1) implies there is some log
canonical centre $Z_{s,t}$ contained in $Z$ that connects $V_s$ and $V_t$.  By symmetry
and (2.2.3) we may find an open subset $U\subset B$ such that $Z_{s,t}$ is a
log canonical centre of $K_Y+\Gamma+2\Gamma_t$.  Varying $s$, $Z_{s,t}$ sweeps out $Z$ and
so by (2.2.2) it follows that $Z$ is a log canonical centre of
$K_Y+\Gamma+2\Gamma_t$.  \qed\enddemo

\proclaim{6.3 Lemma} Let $X$ be a scheme of pure dimension $n$ and let 
$D$ be an ample $\Bbb Q$-divisor such that $L=rD$ is Cartier.  Suppose that there are
effective $\Bbb Q$-divisors $B$ and $R$, where $B\neq 0$ is integral and $D/w\qle aB+R$ 
for some constant $a$.

 Then 
$$
D^n\geq \frac {wa}{r^{n-1}}.
$$
\endproclaim
\demo{Proof} Indeed
$$
r^{n-1}D^n=L^{n-1}\cdot D=w(L^{n-1}\cdot aB+L^{n-1}\cdot R)\geq wa,
$$
where we used the fact that $L^{n-1}\cdot B\geq 1$. \qed\enddemo

\proclaim{6.4 Lemma} Suppose that $K_X+\Delta$ is kawamata log terminal 
of log discrepancy at least $\epsilon>0$, $-(K_X+\Delta)$ is big and nef, where $\Delta$
is an effective $\Bbb Q$-divisor.  Suppose that the dimension of $X$ is $n$ and that
$r(K_X+\Delta)$ is Cartier for some integer $r$.  Assume that the degree is at least
$$
(2n)^n\left (2^nr(n-1)!\right)^{n(n-1)}
$$
and let $\pi\:\map Y.X.$ be the birational morphism determined by a birational family of
tigers of large weight and minimal dimension, whose existence is guaranteed by
(3.3) (here we take $D=-(K_X+\Delta)$).  Define $\Gamma$ to be the log
pullback of $\Delta$ and define $\Gamma_t$ as in (6.1).
\roster 
\item $\Gamma$ is effective in a neighbourhood of the generic fibre of $f$.  
\item There is an open subset $U\subset B$ such that $K_{V_t}+\Gamma_t$ is 
kawamata log terminal of log discrepancy at least $\epsilon$, for all $t\in U$.  
\item If the degree is at least 
$$
\left(\frac {2n}{\epsilon}\right)^n\left (2^nr(n-1)!\right)^{n(n-1)}
$$
then $\pi$ is small in a neighbourhood of the generic fibre of $f$.   
\endroster 
\endproclaim
\demo{Proof} (2) follows from (1) and (6.1).  Thus it suffices to prove 
(1) and (3).  

 Let $\list E.k.$ be the exceptional divisors of $\pi$ which dominate $B$.
By (6.2) $E_i$ has log discrepancy less than or equal to zero with respect to
$K_X+\Delta+2\Delta_t$.  Thus if $a_i$ is the log discrepancy of $E_i$ with
respect to $K_X+\Delta$ then the coefficient of $E_i$ with respect to $\Delta_t$ is at
least $a_i/2$.  By (6.3) the degree of $-(K_{V_t}+\Theta_t)$ is at least
$$
\frac {wa_i}{2r^{n-1}}
$$

 On the other hand by (5.4) the degree of $-(K_{V_t}+\Theta_t)$ is at most $(2^{k+1}k!)^k$. 
Thus 
$$
w\leq \frac {2r^{n-1}(2^{k+1}k!)^k}{a_i}.
$$
By (3.3.2) it follows that 
$$
d\leq (wn)^n\leq \left(\frac {2nr^{n-1}(2^{k+1}k!)^k}{a_i}\right )^n\leq \left(\frac {2n}{a_i}\right)^n\left (2^nr(n-1)!\right)^{n(n-1)}.
$$
Note that if $\Gamma$ is not effective then there exists $i$ such that $a_i\geq 1$ and for
all $i$ we have $a_i>\epsilon$.  (1) and (3) follow immediately.  \qed\enddemo

\demo{Proof of (1.5)} (1.5.1) follows from (3.3.2),
(4.2) and (3.4).  (1.5.2) follows from
(4.3) and (5.4).  (1.5.3) and (1.5.4)
follow from (6.4). \qed\enddemo

\demo{Proof of (1.9)} Immediate from (1.5.1). \qed\enddemo

\demo{Proof of (1.6)} Suppose not.  Then by (1.5) there is a 
birational morphism $\pi\:\map Y.X.$ and a contraction morphism $f\:\map Y.B.$ such that $B$
and the generic fibre have positive dimension.  Moreover $\pi$ is small in a neighbourhood
of the generic fibre of $f$.  Let $H$ be the pullback of an ample divisor on $B$ and let
$G$ be the pushforward of $H$.  Then $G$ is $\Bbb Q$-Cartier as $X$ is $\Bbb Q$-factorial
and so $\pi^*G=H+E$ where $E$ is exceptional and does not dominate $B$.  Thus $G$ is not
big, which contradicts the fact that $X$ has Picard number one.  \qed\enddemo

\proclaim{6.5 Lemma} Let $\pi\:\map Y.B.$ be a bounded family, which parametrises 
triples $(X,L,E)$, where $X$ is a projective variety, $L$ is a line bundle on $X$ and $E$
is a $\Bbb Q$-divisor.  Fix an integer $r$ and a positive real number $\epsilon$.

 Then there is a constant $w_0$ such that for every triple $(X,\Delta,\Gamma)$ where
\roster 
\item $\Delta$ is a subboundary such that $\Delta+E$ is effective. 
\item $K_X+\Delta$ is kawamata log terminal of log discrepancy at least $\epsilon>0$. 
\item $\ring X.(r(K_X+\Delta))=L$. 
\item $K_X+\Delta+\Gamma$ is not kawamata log terminal.   
\item $\Gamma\qle L/w$ 
\endroster 
Then $w\leq w_0$.
\endproclaim
\demo{Proof} Note that divisors $\Delta$ satisfying (1)-(3) form a bounded family.  Thus 
possibly enlarging $B$ we may assume that $B$ parametrises triples $(X,L,\Delta)$.  Note
that we are free to add a base point free line bundle to $L$.  Adding on a sufficiently
ample line bundle to $L$ and passing to a power of $L$ we may therefore assume that $L$ is
very ample.

By Noetherian induction we may assume that $B$ is irreducible and it suffices to exhibit a
non-empty open subset of $B$ with the given property.  Picking a resolution of the
geometric generic fibre and passing to an open subset we may therefore assume that $\pi$
is smooth and that the support of $\Delta$ has global normal crossings.  As we are
assuming that the log discrepancy of $K_X+\Delta$ is at least $\epsilon$ it follows that
every component of $\Delta$ has coefficient at most $1-\epsilon$.

 As $K_X+\Delta+\Gamma$ is not kawamata log terminal (4.1.2) implies that
there is a centre $Z$ of $X$ with multiplicity at least $\epsilon$ in $\Gamma$.  Thus we can find divisors $\list H.k.$ belonging to $|L|$
whose intersection with $Z$ is a finite set of points, where $k$ is the dimension of $Z$.
Now pick a point $p$ of this finite set and choose a further series $H_{k+1},\dots
H_{n-1}$ of elements of $L$ that contain $p$ whose intersection with $\Gamma$ is a finite
set of points.

 Then 
$$
L^n=w(H_1\cdot H_2\cdot \dots H_{n-1}\cdot\Gamma)\geq w\epsilon,
$$
so that 
$$w_0=\frac {L^n}{\epsilon}
$$ 
will do. \qed\enddemo

\demo{Proof of (1.3)} We may as well suppose that the degree of $-(K_X+\Delta)$ is at least
$$
(2n)^n\left (2^nr(n-1)!\right)^{n(n-1)}.
$$
Thus by (1.5) there is a birational morphism $\pi\:\map Y.X.$ and a
contraction morphism $f\:\map Y.B.$ where the fibres of $f$ have degree at most $(2^{k+1}k!)^k$.
A collection of varieties of bounded degree forms a bounded family and so the fibres of
$f$ belong to a bounded family.  Moreover by (6.4) if $\Gamma$ is the log
pullback of $\Delta$ and $\Gamma_t$ the restriction of $\Gamma$ to $V_t$ then we have that
$K_{V_t}+\Gamma_t$ is kawamata log terminal of log discrepancy at least $\epsilon$ for
$t\in U$ an open subset of $B$.

 Set $\Gamma_{t,s}=\pi^*\Delta_s|_{V_t}$.  Then by (6.1.3)
$K_{V_t}+\Gamma_t+\Gamma_{t,s}$ is not kawamata log terminal, for generic choice of $s$
and $t$.  By (6.5) it follows that there is a constant $w_0$ such that $w\leq
w_0$.  Thus by (3.3) the degree of $-(K_X+\Delta)$ is at most $(w_0n)^n$.
\qed\enddemo

\demo{Proof of (1.2)} By the main Theorem of \cite{Kollar93b} given $n$ and $r$ 
there is a fixed integer $M$ such that $L=-M(K_X+\Delta)$ is very ample.  Thus to prove
(1.2) it suffices to bound the degree of $-(K_X+\Delta)$ and we may apply 
(1.3). \qed\enddemo

\demo{Proof of (1.4)} Suppose that $K_X+\Delta$ is kawamata log terminal and 
that $-(K_X+\Delta)$ is big and nef.  Let $r$ be any positive integer such that
$r(K_X+\Delta)$ is Cartier.  By (1.3) there is a constant $M$ such that $d<M$
for all pairs $(X,\Delta)$ satisfying (1.3.1-4).

 Let $U$ be any open subset of $X$ whose complement has codimension at least two.  Suppose
that the algebraic fundamental group of $U$ has cardinality at least $k$.  Then there is a
finite cover $\pi\:\map Y.X.$ of $X$ \'etale over $U$ of degree at least $k$.  Let
$\Gamma$ be the log pullback of $\Delta$.  As $\pi$ is \'etale in codimension one it
follows that $\Gamma$ is effective and the pair $(Y,\Gamma)$ satisfies
(1.3.1-4).  Moreover the degree of $-(K_Y+\Gamma)$ is equal to $k^n$ times the
degree of $-(K_X+\Delta)$.  Thus $k$ is no more than the $n$th root of $M$ divided by the
degree of $-(K_X+\Delta)$.  \qed\enddemo

\head Bibliography \endhead 
\bibliography{/home/mckernan/Tex/math}

\end